\documentclass{amsart}

\usepackage{amssymb,amsmath,leqno}
\usepackage[colorlinks,plainpages,backref]{hyperref}
\usepackage{mathrsfs}
\usepackage[neveradjust]{paralist}
\usepackage[all]{xypic}
\xyoption{dvips}
\usepackage{epic,eepic}
\usepackage{url}

\theoremstyle{definition}
\newtheorem{ntn}{Notation}[section]
\newtheorem{dfn}[ntn]{Definition}
\theoremstyle{plain}

\newtheorem{cnj}[ntn]{Conjecture}

\newtheorem{que}[ntn]{Question}
\theoremstyle{remark}
\newtheorem{rmk}[ntn]{Remark}

\newcommand{\boldone}{{\mathbf{1}}}

\newcommand{\del}{\partial}

\newcommand{\de}{{\mathrm d}}

\newcommand{\eps}{\varepsilon}

\newcommand{\ideal}[1]{{\langle#1\rangle}}
\newcommand{\into}{\hookrightarrow}

\newcommand{\ydx}{y\, {\mathrm d}x}

\newcommand{\calA}{\mathscr{A}}
\newcommand{\calB}{\mathscr{B}}
\newcommand{\calD}{\mathscr{D}}

\newcommand{\calH}{\mathscr{H}}
\newcommand{\calI}{\mathscr{I}}
\newcommand{\calL}{\mathscr{L}}
\newcommand{\calM}{\mathscr{M}}
\newcommand{\calN}{\mathscr{N}}
\newcommand{\calO}{\mathscr{O}}

\newcommand{\fraka}{{\mathfrak{a}}}
\newcommand{\frakm}{{\mathfrak{m}}}

\newcommand{\CC}{\mathbb{C}}
\newcommand{\FF}{\mathbb{F}}
\newcommand{\NN}{\mathbb{N}}

\newcommand{\QQ}{\mathbb{Q}}
\newcommand{\RR}{\mathbb{R}}
\newcommand{\ZZ}{\mathbb{Z}}

\DeclareMathOperator{\ann}{ann}

\DeclareMathOperator{\Der}{Der}
\DeclareMathOperator{\charCycle}{charC}
\DeclareMathOperator{\charVar}{charV}
\DeclareMathOperator{\ev}{ev}

\DeclareMathOperator{\gr}{gr}

\DeclareMathOperator{\ini}{in}

\DeclareMathOperator{\inc}{inc}
\DeclareMathOperator{\Jac}{Jac}

\DeclareMathOperator{\lcm}{lcm}
\DeclareMathOperator{\lct}{lct}

\DeclareMathOperator{\OS}{OS}

\DeclareMathOperator{\Sp}{Sp}

\DeclareMathOperator{\Var}{Var}

\def\mustata{Musta\c t\u a}

\def\BS{Bernstein--Sato}
\def\Grob{Gr\"obner}

\def\GBs{Gr\"obner bases}

\newcommand{\p}{{\partial}}
\def \<{\langle}
\def \>{\rangle}
\newcommand{\monoid}[1]{\left[#1\right]}
\newcommand{\Span}[2]{{#1}\cdot\{#2\}}

\begin{document}

\title[On the $D$-module $f^s$]
{Survey on the $D$-module $f^s$}


\author{Uli Walther\\~\\with an appendix by Anton Leykin}
\address{U. Walther\\ Purdue University\\ Dept.\ of Mathematics\\ 150
  N.\ University St.\\ West Lafayette, IN 47907\\ USA}
\email{walther@math.purdue.edu}

\thanks{This material is based in part upon work supported by the
  National Science Foundation under Grant No.~0932078 000, while UW
  was in residence at the Mathematical Science Research Institute
  (MSRI) in Berkeley, California, during the spring of 2013.}
\thanks{UW was supported in part by the NSF under grants DMS-0901123
  and DMS-1401392}

\address{A. Leykin\\ School of Mathematics\\ Georgia Institute of Technology\\ 686 Cherry St.\\ Atlanta, GA 30332\\ USA}
\email{leykin@math.gatech.edu} 

\thanks{AL was supported in part by the NSF under grant DMS~1151297.}





\maketitle
\setcounter{tocdepth}{1} 
\tableofcontents
\numberwithin{equation}{section}

In this survey we discuss various aspects of the singularity invariants with
differential origin derived from the $D$-module generated
by $f^s$. We should like to point the reader to some other works:
\cite{Saito-Compos07} for $V$-filtration, Bernstein--Sato polynomials,
multiplier ideals; \cite{Budur-survey} for all these and Milnor fibers;
\cite{Torrelli-survey} and \cite{Narvaez-Contemp08} for homogeneity
and free divisors; \cite{Suciu-1301.4851} on details of arrangements,
specifically their Milnor fibers, although less focused on $D$-modules.

We are greatly indebted to Nero Budur, Francisco
J.\ Castro-Jim{\'e}nez, Luis Narv{\'a}ez-Macarro, Morihiko Saito, Wim
Veys and an unknown referee for their careful reading of early
versions of this article and their relentless hunt for errors. We
claim ownership, and apologize for, all surviving mistakes, oversights
and omissions.

\section{Introduction}

\begin{ntn}
In this article, $X$ will denote a complex manifold. Unless indicated
otherwise, $X$ will be $\CC^n$.

Throughout, let $R=\CC[x_1,\ldots,x_n]$ be the ring of polynomials in
$n$ variables over the complex numbers. We denote by $D=R\langle
\del_1,\ldots,\del_n\rangle$ the Weyl algebra. In particular, $\del_i$
denotes the partial differentiation operator with respect to $x_i$. If
$X$ is a general manifold, $\calO_X$ (the sheaf of regular functions)
and $\calD_X$ (the sheaf of $\CC$-linear differential operators on
$\calO_X$) take the places of $R$ and $D$.

If $X=\CC^n$ we use Roman letters to denote rings and modules; in the
general case we use calligraphic letters to denote corresponding sheaves.

By the ideal $J_f$ we mean the $\calO_X$-ideal generated by
the partial derivatives $\frac{\del f}{\del x_1},\ldots,\frac{\del
  f}{\del x_n}$;
this ideal varies with the choice of coordinate system in which we
calculate. In contrast, the Jacobian ideal $\Jac(f)=J_f+(f)$ is
independent.
\end{ntn}

The ring $D$ (resp.~the sheaf $\calD_X$) is coherent, and both left-
and right-Noetherian; it has only trivial two-sided ideals
\cite[Thm.~1.2.5]{Bjoerk-anal}. Introductions to the theory of
$D$-modules as we use them here can be found in
\cite{Kashiwara-AMSbook,Bernstein-notes,Bjoerk-anal,Bjoerk-rings}.

The ring $D$ admits the order filtration induced by the weight $x_i\to
0$, $\del_i\to 1$.  The order filtration (and other good filtrations)
leads to graded objects $\gr_{(0,1)}(-)$, see \cite{Schapira85}. The
graded objects obtained from ideals are ideals in the polynomial ring
$\CC[x,\xi]$, homogeneous in the symbols of the differentiation
operators; their radicals are closed under the Poisson bracket, and
thus the corresponding varieties are involutive
\cite{Kashiwara-RIMS-I,KashiwaraKawai-RIMS-III}. For a $D$-module $M$
and a component $C$ of the support of $\gr_{(0,1)}(M)$, attach to the
pair $(M,C)$ the multiplicity $\mu(M,C)$ of $\gr_{(0,1)}(M)$ along
$C$. The characteristic cycle of $M$ is $\charCycle(M)=\sum_C
\mu(M,C)\cdot C$, an element of the Chow ring on $T^*\CC^{n}$. The
module is \emph{holonomic} if it is zero or if its characteristic
variety is of dimension $n$, the minimal possible value.

Throughout, $f$ will be a regular function on $X$, with divisor
$\Var(f)$. We distinguish several homogeneity conditions on $f$:
\begin{itemize}
\item $f$ is \emph{locally} (\emph{strongly}) \emph{Euler-homogeneous}
  if for all $p\in\Var(f)$ there is a vector field $\theta_p$ defined
  near $p$ with $\theta_p\bullet(f)=f$ (and $\theta_p$ vanishes at
  $p$).
\item $f$ is \emph{locally }(\emph{weakly}) \emph{quasi-homogeneous}
  if near all $p\in \Var(f)$ there is a local coordinate system
  $\{x_i\}$ and a positive (resp.~non-negative) weight vector
  $a=\{a_1,\ldots,a_n\}$ with respect to which $f=\sum_{i=1}^n a_i
  x_i\del_i(f)$.
\item We reserve \emph{homogeneous} and
  \emph{quasi-homogeneous} for the case when $X=\CC^n$ and $f$ is
  globally homogeneous or quasi-homogeneous.
\end{itemize}

To any non-constant $f\in R$, one can attach
several invariants that measure the singularity structure of the
hypersurface $f=0$. In this article, we are primarily interested in
those derived from the (parametric) annihilator $\ann_{D[s]}(f^s)$ of $f^s$:
\begin{dfn}
 Let $s$ be a new variable, and denote by $R_f[s]\cdot f^s$ the free module
 generated by $f^s$ over the localized ring $R_f[s]=R[f^{-1},s]$. Via the
 chain rule
\begin{eqnarray}\label{eq-f^s}
\del_i\bullet
(\frac{g}{f^k}f^s)=\del_i\bullet(\frac{g}{f^k})f^s+
\frac{sg}{f^{k+1}}\cdot\frac{\del f}{\del x_i}f^s
\end{eqnarray}
for each $g(x,s)\in R[s]$, $R_f[s]\cdot f^s$ acquires the structure of a
left $D[s]$-module. Denote by
\[
\ann_{D[s]}(f^s)=\{P\in D[s]\mid P\bullet f^s=0\}
\]
the \emph{parametric
  annihilator}, and by
\[
\calM_f(s)=D[s]/\ann_{D[s]}(f^s)
\]
the cyclic $D[s]$-module generated by $1\cdot f^s\in
R_f[s]\cdot f^s$.
\end{dfn}

\medskip

Bernstein's functional equation \cite{Bernstein-Nauk72}
asserts the existence of a differential operator
$P(x,\del,s)$ and a nonzero polynomial $b_{f,P}(s)\in\CC[s]$ such that
\begin{eqnarray}\label{eqn-bernstein}
P(x,\del,s)\bullet f^{s+1}=b_{f,P}(s)\cdot f^s,
\end{eqnarray}
\emph{i.e.}\ the existence of the element $P\cdot
f-b_{f,P}(s)\in\ann_{D[s]}(f^s)$. 
Bernstein's result implies that $D[s]\bullet f^s$ is $D$-coherent
(while $R_f[s]f^s$ is not).
\begin{dfn}
The monic generator of the ideal in $\CC[s]$ generated by all $b_{f,P}(s)$
appearing in an equation \eqref{eqn-bernstein} is the
\emph{Bernstein--Sato polynomial} $b_f(s)$.
Denote $\rho_f\subseteq \CC$ the set of roots of $b_f(s)$.
\end{dfn}
Note that the operator $P$ in the functional equation is only
determined up to $\ann_{D[s]}(f^s)$.  See \cite{Bjoerk-rings} for an
elementary proof of the existence of $b_f(s)$. Alternative (and more
general) proofs are given in \cite{Kashiwara-AMSbook}; see also
\cite{Bernstein-notes,MebkhoutNavarro-Ecole91,Nunez-JA13}.

The $\CC[s]$-module $\calM_f(s)/\calM_f(s+1)$ is precisely annihilated
by $b_f(s)$. It is an interesting problem to determine for any
$q(s)\in\CC[s]$ the ideals $\fraka_{f,q(s)}=\{g\in R\mid q(s)g f^s\in
D[s]\bullet f^{s+1}\}$ from \cite{Walther-Bernstein}.
By \cite{Malgrange-isolee},
$\fraka_{f,s+1}=R\cap(\ann_{D[s]}(f^s)+D[s]\cdot (f,J_f))$.
\begin{que}
Is $\fraka_{f,s+1}=J_f +(f)$?
\end{que}
A positive answer would throw light on connections between $b_f(s)$
and cohomology of Milnor fibers.
\begin{rmk}
At the 1954 International Congress of Mathematics in Amsterdam,
I.M.~Gel'fand asked the following question. Given a real analytic
function $f\colon \RR^n\to\RR$, the assignment ($s\in\CC$)
\[
f(x)^s_+=\left\{\begin{array}{ll}f(x)^s&\text{ if }f(x)>0,\\0&\text{
  if }f(x)\le 0\end{array}\right.
\]
is continuous in $x$ and analytic in $s$ where the real part of $s$ is
positive. Can one analytically continue $f(x)^s_+$? Sato introduced
$b_f(s)$ in order to answer Gel'fand's question; Bernstein
\cite{Bernstein-Nauk72} established their existence in general.
\end{rmk}
\begin{rmk}
Let $m\in M$ be a nonzero section of a holonomic
$D$-module. Generalizing the case $1\in R$ there is a functional
equation
\[
P(x,\del,s)\bullet (mf^{s+1})=b_{f,P;m}(s)\cdot mf^s
\]
with $b_{f,P;m}(s)\in\CC[s]$ nonzero. The monic generator of
the ideal $\{b_{f,P;m}(s)\}$ is the $b$-function $b_{f;m}(s)$,
\cite{Kashiwara-bfu}.
\end{rmk}

\section{Parameters and numbers}

For any complex number $\gamma$, the expression $f^\gamma$ represents,
locally outside $\Var(f)$, a multi-valued analytic function. Via the
chain rule as in \eqref{eq-f^s}, the cyclic $R_f$-module $R_f\cdot
f^\gamma$ becomes a left $D$-module, and we set
\[
\calM_f(\gamma)=D\bullet f^\lambda\cong D/\ann_D(f^\gamma).
\]
There are natural
$D[s]$-linear maps
\[
\ev_f(\gamma)\colon \calM_f(s)\to \calM_f(\gamma),\qquad
P(x,\del,s)\bullet f^s\mapsto P(x,\del,\gamma)\bullet f^\gamma,
\]
and $D$-linear inclusions
\[
\inc_f(s)\colon \calM_f(s+1)\to \calM_f(s),\qquad
P(x,\del,s)\bullet f^{s+1}\mapsto P(x,\del,s)\cdot f\bullet f^s
\]
with cokernel $\calN_f(s)=\calM_f(s)/\calM_f(s+1)\cong
D[s]/(\ann_{D[s]}(f^s)+D[s]f)$, and
\[
\inc_f(\gamma)\colon \calM_f(\gamma+1)\to \calM_f(\gamma),\qquad
P(x,\del)\bullet f^{\lambda+1}\mapsto P(x,\del)\cdot f\bullet f^\lambda
\]
with cokernel $\calN_f(\gamma)=\calM_f(\gamma)/
\calM_f(\gamma+1)\cong D/(\ann_D(f^\gamma)+D\cdot f)$.

The kernel of the morphism $\ev_f(\gamma)$ contains the (two-sided)
ideal $D[s](s-\gamma)$; the containment can be proper, for example if
$\gamma=0$.  If $\{ \gamma-1,\gamma-2,\ldots\}$ is disjoint from the
root set $\rho_f$ then $\ker \ev_f(\gamma)=D[s]\cdot(s-\gamma)$,
\cite{Kashiwara-bfu}. If $\gamma\not\in\rho_f$ then $\inc_f(\gamma)$
is an isomorphism because of the functional equation; if $\gamma=-1$,
or if $b_f(\gamma)=0$ while $\rho_f$ does not meet
$\{\gamma-1,\gamma-2,\ldots\}$ then $\inc_f(\gamma)$ is not surjective
\cite{Walther-Bernstein}.
\begin{que}
Does $\inc_f(\gamma)$ fail to be an isomorphism for all
$\gamma\in\rho_f$?\footnote{In a recent preprint, this question is
  answered in the negative by M.~Saito, see \cite{Saito-powers}.}
\end{que}
In contrast, the induced maps $\calM_f(s)/(s-\gamma-1)\to
\calM_f(s)/(s-\gamma)$ are isomorphisms exactly when $\gamma\not\in
\rho_f$, \cite[6.3.15]{Bjoerk-anal}.
The morphism $\inc_f(s)$ is never surjective as $s+1$ divides
$b_f(s)$. One sets
\[
\tilde b_f(s)=\frac{b_f(s)}{s+1}.
\]

By \cite[4.2]{Torrelli-RIMS09}, the following are equivalent for
a section $m\not =0$ of a holonomic module:
\begin{itemize}
\item the smallest integral root of $b_{f;m}(s)$ is at least $-\ell$;
\item $(D\bullet m)\otimes_RR[f^{-1}]$ is generated by $m/f^\ell
=m\otimes 1/f^\ell$;
\item $(D\bullet m)\otimes_RR[f^{-1}]/D\bullet (m\otimes 1)$ is
  generated by $m/f^\ell$;
\item $D[s]\bullet mf^s\to (D\bullet m)\otimes_RR[f^{-1}]$, $P(s)\bullet
  (mf^s)\mapsto P(-\ell)\bullet (m/f^\ell)$ is an epimorphism with
  kernel $D[s]\cdot(s+\ell)mf^s$.
\end{itemize}

\begin{dfn}
We say that  $f$ satisfies condition
\begin{itemize}
\item $(A_1)$ (resp.~$(A_s)$) if $\ann_D(1/f)$ (resp.\ $\ann_D(f^s)$)
  is generated by operators of order one;
\item $(B_1)$ if $R_f$ is generated by $1/f$ over $D$.
\end{itemize}
\end{dfn}

Condition $(A_1)$ implies $(B_1)$ in any case \cite{Torrelli-Bull04}.
Local Euler-homogeneity, $(A_s)$ and $(B_1)$ combined imply $(A_1)$
\cite{Torrelli-survey}, and for Koszul free divisors (see
Definition~\ref{dfn-Koszul-free} below) this implication can be
reversed \cite{Torrelli-Bull04}.

Condition $(A_1)$ does not imply $(A_s)$: $f=xy(x+y)(x+yz)$ is free
(see Definition \ref{dfn-free}), and locally Euler-homogeneous and satisfies
$(A_1)$ and $(B_1)$
\cite{Calderon-Ecole99,CalderonCastroMondNarvaez-Helvetici02,
  CalderonNarvaez-Compos02,CastroUcha-JSC01,Torrelli-Bull04}, but
$\ann_{D[s]}(f^s)$ and $\ann_D(f^s)$ require a second order generator.

Condition $(A_1)$ implies local Euler-homogeneity if $f$ has isolated
singularities \cite{Torrelli-Fourier02}, or if it is Koszul-free or of
the form $z^n-g(x,y)$ for reduced $g$ \cite{Torrelli-Bull04}. In
\cite{CastroGagoHartilloUcha-Revista07} it is shown that for certain
locally weakly quasi-homogeneous free divisors $\Var(f)$, $(A_1)$
holds for high powers of $f$, and even for $f$ itself by
\cite[Rem.~1.7.4]{Narvaez-Contemp08}.

For an isolated singularity, $f$ has $(A_1)$ if and only if it has
$(B_1)$ and is quasi-homogeneous \cite{Torrelli-Fourier02}. For
example, a reduced plane curve (has automatically $(B_1)$ and) has
$(A_1)$ if and only if it is quasi-homogeneous. See
\cite{Schulze-PAMS07} for further results.

Condition $(B_1)$ is equivalent to $\inc_f(-2), \inc_f(-3),\ldots $
all being isomorphisms, and also to $-1$ being the only integral root
of $b_f(s)$, \cite{Kashiwara-bfu}.  Locally quasi-homogeneous free
divisors satisfy condition $(B_1)$ at any point,
\cite{CastroUcha-Steklov02}.

\section{$V$-filtration and Bernstein--Sato polynomials}

\subsection{$V$-filtration}

The articles
\cite{Saito-Bull94,MaisonobeMebkhout-Paris04,Budur-Vfilt05,Budur-survey}
are recommended for material on $V$-filtrations.

\subsubsection{Definition and basic properties}
Let $Y$ be a smooth complex manifold (or variety), and let $X$ be a
closed submanifold (or -variety) of $Y$ defined by the ideal sheaf
$\calI$. The $V$-filtration on $\calD_Y$ along $X$ is, for $k\in\ZZ$, given by
\[
V^k(\calD_Y)=\{P\in \calD_Y\mid P\bullet \calI^{k'}\subseteq \calI^{k+k'}\quad\forall
k'\in\ZZ\}
\]
with the understanding that $\calI^{k'}=\calO_Y$ for $k'\le 0$.
The associated graded
sheaf of rings $\gr_V(\calD_Y)$ is isomorphic to the sheaf of rings of
differential operators on the normal bundle $T_X(Y)$, algebraic in the
fiber of the bundle.

Suppose that $Y=\CC^n\times\CC$ with coordinate function
$t$ on $\CC$, and let $X$ be the hyperplane $t=0$. Then $V^{k}(D_Y)$
is spanned by $\{x^u\del^vt^a\del_t^b\mid a-b\geq k\}$.  Given a
coherent holonomic $D_Y$-module $M$ with regular singularities in the
sense of \cite{KashiwaraKawai-regSing81}, Kashiwara and Malgrange
\cite{Malgrange-evan,Kashiwara-LNM1016} define an exhaustive
decreasing rationally indexed filtration on $M$ that is compatible
with the $V$-filtration on $D_Y$ and has the following properties:
\begin{enumerate}
\item each $V^\alpha(M)$ is coherent over $V^0(D_Y)$ and the set of
  $\alpha$ with nonzero $\gr^\alpha_V(M)=V^\alpha(M)/V^{>\alpha}(M)$
  has no accumulation point;
\item for $\alpha\gg 0$, $V^1(D_Y)V^\alpha(M)=V^{\alpha+1}(M)$;
\item $t\del_t-\alpha$ acts nilpotently on $\gr^\alpha_V(M)$.
\end{enumerate}
The $V$-filtration is unique and can be defined in somewhat greater
generality \cite{Budur-Vfilt05}. Of special interest is the following
case considered in \cite{Malgrange-evan,Kashiwara-LNM1016}.
\begin{ntn}
Denote $R_{x,t}$ the polynomial ring $R[t]$, $t$ a new indeterminate,
  and let $D_{x,t}$ be the corresponding Weyl algebra.
Fix $f\in R$ and consider the regular $D_{x,t}$-module
\[
\calB_f=H^1_{f-t}(R[t]),
\]
the unique local cohomology module of $R[t]$ supported in $f-t$. Then
$\calB_f$ is naturally isomorphic as $D_{x,t}$-module to the direct
image (in the $D$-category)  $i_+(R)$ of $R$  under the graph embedding
\[
i\colon X\to X\times \CC,\qquad
x\mapsto (x,f(x)).
\]
Moreover, extending \eqref{eq-f^s} via
\[
t\bullet (g(x,s)f^{s-k})=g(x,s+1)f^{s+1-k};\qquad
\del_t\bullet (g(x,s)f^{s-k})=-sg(x,s-1)f^{s-1-k},
\]
the module $R_f[s]\otimes f^s$ becomes a $D_{x,t}$-module extending the
$D[s]$-action  where $-\del_tt$
acts as $s$.
\end{ntn}

The existence of the $V$-filtration on
$\calB_f=i_+(R)$ is equivalent to the existence of generalized $b$-functions
$b_{f;\eta}(s)$ in the sense of \cite{Kashiwara-bfu}, see
\cite{Kashiwara-holonomicII,Malgrange-evan}.  In fact, one can recover
one from the other:
\[
V^\alpha(\calB_f)=\{\eta\in \calB_f\mid [b_{f;\eta}(-c)=0]\Rightarrow
[\alpha\le c]\}
\]
and the multiplicity of $b_{f;\eta}(s)$ at $\alpha$ is the degree of
the minimal polynomial of $s-\alpha$ on $\gr_V^\alpha(D[s]\eta
f^s/D[s]\eta f^{s+1})$, \cite{Sabbah-Rabida87}. For more on this
``microlocal approach'' see \cite{Saito-Bull94}.

\subsection{The log-canonical threshold}
By \cite{Kollar-Proc97}, see also \cite{Lichtin-Arkiv89,Yano}, the
absolute value of the largest root of $b_f(s)$ is the
\emph{log-canonical threshold} $\lct(f)$ given by the supremum of all
numbers $s$ such that the local integrals
\[
\int_{U\ni p}\frac{|dx|}{|f|^{2s}}
\]
converge for all $p\in X$ and all small open $U$ around $p$.  Smaller
$\lct$ corresponds to worse singularities; the best one
can hope for is $\lct(f)=1$ as one sees by looking at a smooth
point. The notion goes back to Arnol'd, who called it (essentially) the
complex singular index \cite{ArnoldGuseinVarchenko-SinI-II85}.

The point of \emph{multiplier ideals} is to force the
finiteness of the integral by allowing moderating functions in the
integral:
\[
\calI(f,\lambda)_p=\{g\in \calO_X\mid \frac{g}{f^\lambda} \text{ is
  $L^2$-integrable near $p\in\Var(f)$}\}
\]
for $\lambda\in\RR$.
By \cite{EinLazarsfeldSmithVarolin-Duke04}, there is a finite
collection of \emph{jumping numbers} for $f$
of rational numbers $0=\alpha_0<\alpha_1<\cdots <\alpha_\ell=1$ such that
$\calI(f,\alpha)$ is constant on $[\alpha_i,\alpha_{i+1})$ but
  $\calI(f,\alpha_i)\not =\calI(f,\alpha_{i+1})$. The log-canonical
  threshold appears as $\alpha_1$. These ideas had
  appeared previously in \cite{Lipman-Dekker82,LoeserVaquie-Top90}.

Generalizing Kollar's approach, each $\alpha_i$ is a root of $b_f(s)$,
\cite{EinLazarsfeldSmithVarolin-Duke04}. In \cite[Thm.~4.4]{Saito-Compos07} a
partial converse is shown for locally Euler-homogeneous divisors.
Extending the idea of jumping numbers to the range $\alpha>1$ one sees
that $\alpha$ is a jumping number if and only if $\alpha+1$ is a
jumping number, but the connection to the Bernstein--Sato polynomial
is lost in general. For example, if $f(x,y)=x^2+y^3$ then jumping
numbers are $\{5/6,1\}+\NN$ while $b_f(s)=(s+5/6)(s+1)(s+7/6)$.

\subsection{Bernstein--Sato polynomial}\label{subsec-bfu}

The roots of $b_f(s)$ relate to an astounding number of other
invariants, see for example \cite{Kollar-Proc97} for a
survey. However, besides the functional equation there is no known way
to describe $\rho_f$.

\subsubsection{Fundamental results}
Let $p\in\CC^n$ be a closed point, cut out by the maximal ideal
$\frakm\subseteq R$. Extending $R$ to the localization $R_\frakm$ (or even the
ring of holomorphic functions at $p$) one arrives at potentially
larger sets of polynomials $b_{f,P}(s)$ that satisfy a functional
equation \eqref{eqn-bernstein} with $P(x,\del,s)$ now in the
correspondingly larger ring of differential operators. The
\emph{local} (resp.\ \emph{local analytic}) Bernstein--Sato polynomial
$b_{f,p}(s)$ (resp.\ $b_{f,p^{an}}(s)$) is the generator of the
resulting ideal generated by the $b_{f,P}(s)$ in $\CC[s]$. We denote by
$\rho_{f,p}$ (resp.\ $\tilde\rho_{f,p}$) the root set of
$b_{f,p^{an}}(s)$ (resp.\ $b_{f,p^{an}}(s)/(s+1)$).
From the
definitions and
\cite{Gennady-PAMS97,BrianconMaisonobe-Tokyo00,BrianconMaynadier-Kyoto99}
\begin{gather}\label{eq-local-bfu}
b_{f,p^{an}}(s)\big| b_{f,p}(s)\big|
b_f(s)=\lcm_{p\in\Var(f)}b_{f,p}(s)=\lcm_{p\in\Var(f)}b_{f,p^{an}}(s),
\end{gather}
and the function $\CC^n\ni p\mapsto \Var(b_f(s))$, counting with
multiplicity, is upper semi-continuous in the sense that for $p'$
sufficiently near $p$ one has $b_{f,p'}(s)|b_{f,p}(s)$.  The
underlying reason is the coherence of $D$.

The Bernstein--Sato polynomial $b_f(s)$ factors over $\QQ$ into linear
factors, $\rho_f\subseteq \QQ$, and all roots are negative
\cite{Malgrange-isolee,Kashiwara-bfu}. The proof uses
resolution of singularities over $\CC$ in order to reduce to simple
normal crossing divisors, where rationality and negativity of the
roots is evident. For this Kashiwara proves a comparison theorem
\cite[Thm.~5.1]{Kashiwara-bfu}  that establishes $b_f(s)$ as a
divisor of a shifted product of the least common multiple of the local
Bernstein--Sato polynomials of the pullback of $f$ under the
resolution map. There is a  refinement by Lichtin
\cite{Lichtin-Arkiv89} for plane curves.  The roots of $b_f(s)$,
besides being negative, are always greater than $-n$, $n$ being the
minimum number of variables required to express $f$ locally
analytically \cite{Varchenko-Izvestiya81,Saito-Bull94}.

\subsubsection{Constructible sheaves from $f^s$}\label{subsec-strat-bfu}
Let $V=V(n,d)$ be the vector space of all complex polynomials in
$x_1\ldots,x_n$ of degree at most $d$. Consider the function
$\beta\colon V\ni f\mapsto b_f(s)$. By
\cite{Gennady-PAMS97,BrianconMaynadier-Kyoto99}, there is an algebraic
stratification of $V$ such that on each stratum the function $\beta$
is constant. For varying $n,d$ these stratifications can be made to be
compatible.

\subsubsection{Special cases}
\label{subsubsec-bfu-examples}
If $p$ is a smooth point of $\Var(f)$ then $f$ can be used as an
analytic coordinate near $p$, hence $b_{f,p^{an}}(s)=s+1$, and so
$b_f(s)=s+1$ for all smooth hypersurfaces. By Proposition 2.6 in
\cite{BrianconMaisonobe-Progr134}, an extension of
\cite{BrianconLaurentMaisonobe-Comptes91}, the equation $b_f(s)=s+1$
implies smoothness of $\Var(f)$. Explicit formul\ae\ for the
Bernstein--Sato polynomial are rare; here are some classes of
examples.
\begin{itemize}
\item $f=\prod x_i^{a_i}$:\quad $P=\prod \del_i^{a_i}$ up to a scalar,
  $b_f(s)=\prod_i\prod_{j=1}^{a_i}(s+j/a_j)$.
\item $f$ (quasi-)homogeneous with isolated singularity at zero:
  $\tilde b_f(s)=\lcm(s+\frac{\deg(g\de x)}{\deg(f)}$,
  where $g$ runs through a (quasi-)homogeneous standard basis for $J_f$
  by work of 
  Kashiwara, Sato, Miwa,
  Malgrange, Kochman \cite{Malgrange-isolee,Yano,Torrelli-JA05,
    Kochman-PNAS76}. Note that the Jacobian ring of such a singularity
  is an Artinian Gorenstein ring, whose duality operator implies
  symmetry of $\rho_f$.
\item $f=\det(x_{i,j})_1^n$:\quad $P=\det(\del_{i,j})_1^n$,
  $b_f(s)=(s+1)\cdots (s+n)$. This is  attributed to Cayley, but see
  the comments in \cite{CarracioloSokalSportiello-AdvAM13}.
\item For some hyperplane arrangements, $b_f(s)$ is known, see
  \cite{Walther-Bernstein,BudurSaitoYuzvinsky-JLMS11}.
\item There is a huge list of examples worked out in \cite{Yano}.
\end{itemize}
If $V$ is a complex vector space, $G$ a reductive group acting
linearly on $V$ with open orbit $U$ such that $V\smallsetminus U$ is a
divisor $\Var(f)$, Sato's theory of prehomogeneous vectors spaces
\cite{SatoShintani-Annals74,Muro-Academic88,SatoShintaniMuro-Nagoya90,Yano-RIMS76}
yields a factorization for $b_f(s)$.  For reductive linear free
divisors, \cite{GrangerSchulze-RIMS10,Sevenheck-Fourier11} discuss
symmetry properties of Bernstein--Sato polynomials.  In
\cite{Narvaez-1201.3594} this theme is taken up again, investigating
specifically symmetry properties of $\rho_f$ when $D[s]\bullet f^s$
has a Spencer logarithmic resolution (see \cite{CastroUcha-Steklov02}
for definitions). This covers locally quasi-homogeneous free divisors,
and more generally free divisors whose Jacobian is of linear type. The
motivation is the fact that roots of $b_f(s)$ seem to come in strands,
and whenever roots can be understood the strands appear to be linked
to Hodge-theory.

There are several results on $\rho_f$ for other divisors of special shape.
Trivially, if $f(x)=g(x_1,\ldots,x_k)\cdot h(x_{k+1},\ldots,x_n)$ then
$b_f(s)\mid b_g(s)\cdot b_h(s)$; the question of equality appears to
be open. In contrast, $b_f(s)$ cannot be assembled
from the Bernstein--Sato polynomials of the factors of $f$
in general, even if the factors are hyperplanes and one has some
control on the intersection behavior, see Section~\ref{sec-arr} below.  If
$f(x)=g(x_1,\ldots,x_k)+h(x_{k+1},\ldots,x_n)$ and at least one is
locally Euler-homogeneous then there are
Thom--Sebastiani type formul\ae\, \cite{Saito-Bull94}. In particular,
diagonal hypersurfaces are completely understood.

\subsubsection{Relation to intersection homology module}
Suppose $Y=\Var(f_1,\ldots,f_k)\subseteq X$ is a complete intersection
and denote by $\calH^k_Y(\calO_X)$ the unique (algebraic) local
cohomology module of $\calO_X$ along $Y$. Brylinski--Kashiwara
\cite{Brylinski-Asterisque83-1,Brylinski-Asterisque83-2} defined
$\calL(Y,X)\subseteq \calH^k_{Y}(\calO_X)$, the \emph{intersection
  homology $\calD_X$-module of $Y$}, the smallest $\calD_X$-module
equal to $\calH^k_{Y}(\calO_X)$ in the generic point(s). See also
\cite{BarletKashiwara-Invent86}. The module $\calL(X,Y)$ contains the
fundamental class of $Y$ in $X$ \cite{Barlet-LNM807}.
\begin{que}
When is $\calL(X,Y)=\calH^k_{Y}(\calO_X)$?
\end{que}
Equality is equivalent to
$\calH^k_{Y}(\calO_X)$ being generated by the cosets of
$\Delta/\prod_{i=1}^k f_i$ over $\calD_X$ where $\Delta$ is the ideal
generated by the $k$-minors of the Jacobian matrix of
$f_1,\ldots,f_k$. A necessary condition is that $1/\prod_{i=1}^k f_i$
generates $\calH^k_{Y}(\calO_X)$, but this is not sufficient:
consider $xy(x+y)(x+yz)$, where $\rho_f=-\{1/2,3/4,1,1,1,5/4\}$.
Indeed, by \cite{Torrelli-RIMS09},
equality can be characterized in terms of functional equations, as the
following are equivalent at $p\in X$:
\begin{enumerate}
\item $\calL(X,Y)=\calH^k_{Y}(\calO_X)$ in the stalk;
\item $\tilde\rho_{f,p}\cap\ZZ=\emptyset$;
\item $1$ is not an eigenvalue of the monodromy operator on
  the reduced cohomology of the Milnor fibers near $p$.
\end{enumerate}
If $1/\prod_{i=1}^k f_i$ generates $R[1/\prod f_i]$ and
$1/\prod_{i=1}^k f_i\in\calL(X,Y)$ then $\tilde b_f(-1)\not =0$,
\cite{Torrelli-RIMS09}.  It seems unknown whether (irrespective of
$1/\prod_{i=1}^k f_i$ generating $R[1/\prod f_i]$) the condition
$\tilde b_f(-1)\not =0$ is equivalent to $1/\prod_{i=1}^k f_i$ being
in $\calL(X,Y)$. See also \cite{Massey-IJM09} for a topological
viewpoint (by the Riemann--Hilbert correspondence of Kashiwara and
Mebkhout \cite{Kashiwara-RIMS84,Mebkhout-Compos84}, $\calL(X,Y)$
corresponds to the intersection cohomology complex of $Y$ on $X$
\cite{Brylinski-Asterisque83-1} and $\calH^k_{Y}(\calO_X)$ to
$\CC_Y[n-k]$, \cite{Grothendieck-deRham,Kashiwara-bfu,Mebkhout-RIMS76};
equality then says: the link is a rational homology sphere).  In
\cite{Barlet-RIMS99}, Barlet characterizes property (3) above in terms
of currents for complexified real $f$. Equivalence of (1) and (3)
for isolated singularities can be derived from
\cite{Milnor,Brieskorn-Manus70}; the general case can be shown using
\cite[4.5.8]{Saito-RIMS90} and the formalism of weights.  For the case
$k=1$, (1) requires irreducibility; in the general case, there is a
criterion in terms of $b$-functions \cite[1.6, 1.10]{Torrelli-RIMS09}.

\section{LCT and logarithmic ideal}\label{sec-LCT}

\subsection{Logarithmic forms}
Let $X=\CC^n$ be the analytic manifold, $f$ a holomorphic function on
$X$, and $Y=\Var(f)$ a divisor in $X$ with $j\colon U=X\smallsetminus
Y\into X$ the embedding. Let $\Omega^\bullet_X(*Y)$ denote the complex
of differential forms on $X$ that are (at worst) meromorphic along
$Y$. By \cite{Grothendieck-deRham}, $\Omega^\bullet_X(*Y)\to \RR
j_*\CC_U$ is a quasi-isomorphism.

A form $\omega$ is \emph{logarithmic} along $Y$ if $f\omega$ and
$f\de\omega$ are holomorphic; these $\omega$ form the logarithmic de
Rham complex $\Omega^\bullet_X(\log Y)$ on $X$ along $Y$. The complex
$\Omega^\bullet_X(\log Y)$ was first used with great effect on normal
crossing divisors by Deligne \cite{Deligne-HodgeII} in order to
establish mixed Hodge structures, and later by Esnault and Viehweg in
order to prove vanishing theorems \cite{EsnaultViehweg-book}. A major
reason for the success of normal crossings is that in that case
$\Omega^i_X(\log Y)$ is a locally free module over $\calO_X$.  The
logarithmic de Rham complex was introduced in \cite{KSaito-Tokyo} for
general divisors.

\subsection{Free divisors}

\begin{dfn}\label{dfn-free}
A divisor $\Var(f)$ is \emph{free} if (locally) $\Omega^1_X(\log f)$ is a
free $\calO_X$-module.
\end{dfn}
For a non-smooth locally Euler-homogeneous divisor, freeness is
equivalent to the Jacobian ring $\calO_X/J_f$ being a Cohen--Macaulay
$\calO_X$-module of codimension $2$; in general, freeness is equivalent
to the Tjurina algebra $R/(f,\frac{\del f}{\del x_1},\ldots,\frac{\del
  f}{\del x_n})$ being of projective dimension $2$ or less over $R$.
See \cite{KSaito-Tokyo, Aleksandrov-Euler86} for relations to
determinantal equations. Free divisors have rather big singular locus,
and are in some ways at the opposite end 
from
isolated singularities
in the singularity zoo.  If $\Omega^1_X(\log f)$ is (locally) free,
then $\Omega^i_X(\log f)\cong \bigwedge^i\Omega^i_X(\log f)$ and also
(locally) free, \cite{KSaito-Tokyo}. A weakening is
\begin{dfn}
A divisor $\Var(f)$ is \emph{tame} if, for all $i\in\NN$, (locally)
$\Omega^i_X(\log f)$ has projective dimension at most $i$ as a
$\calO_X$-module.
\end{dfn}

Plane curves are
trivially free; surfaces in $3$-space are trivially tame.  Normal
crossing divisors are easily shown to be free.  Discriminants of
(semi)versal deformations of an isolated complete intersection
singularity  (and some others) are free,
\cite{Aleksandrov-Euler86,Aleksandrov-Soviet90,Looijenga-Isol84, KSaito-unfolding,Damon-AJM98,BuchweitzEbelingVonBothmer-JAG09}.
Unitary reflection arrangements are free \cite{Terao-Invent81}.

\begin{dfn}
The \emph{logarithmic derivations} $\Der_X(-\log f)$ along $Y=\Var(f)$
are the $\CC$-linear derivations $\theta\in\Der(\calO_X;\CC)$ that
satisfy $\theta\bullet f\in(f)$.
\end{dfn}
A derivation $\theta$ is logarithmic along $Y$ if and only it is so
along each component of the reduced divisor to $Y$ \cite{KSaito-Tokyo}. The
modules $\Der_X(-\log f)$ and $\Omega^1_X(\log f)$ are reflexive and
mutually dual over $R$. Moreover, $\Omega^i_X(\log f)$ and
$\Omega^{n-i}_X(\log f)$ are dual.

\subsection{LCT}
\begin{dfn}
If
\begin{eqnarray}\label{eqn-LCT}
\Omega^\bullet_X(\log Y)\to \Omega^\bullet_X(*Y)
\end{eqnarray}
is a quasi-isomorphism, we say that \emph{LCT holds for $Y$}.
\end{dfn}
We recommend \cite{Narvaez-Contemp08}.

\begin{rmk}
\begin{asparaenum}
\item This ``Logarithmic Comparison Theorem'', a property of a
  divisor, is very hard to check explicitly. No general algorithms are
  known, even in $\CC^3$ (but see \cite{CastroTakayama-JAlg09} for $n=2$).
\item LCT fails for rather simple divisors such as
  $f=x_1x_2+x_3x_4$.
\item If $Y$ is a reduced normal crossing divisor, Deligne proved
  \eqref{eqn-LCT} to be a filtered (by pole filtration) quasi-isomorphism
  \cite{Deligne-LNM163}; this provided a crucial step in the
  development of the theory of mixed Hodge structures
  \cite{Deligne-HodgeII}.
\item Limiting the order of poles in forms needed to capture all
  cohomology of $U$ started with the seminal article
  \cite{Griffiths-Annals69} and continues, see for example
  \cite{DeligneDimca-Ecole90,Dimca-AJM91,Karpishpan-Compos91}.
\item The free case was studied for example in
  \cite{CastroMondNarvaez-TAMS96}. But even in this case, LCT is not
  understood.
\item If $f$ is quasi-homogeneous with an isolated singularity at
  the origin, then LCT for $f$ is equivalent to a topological
  condition (the link of $f$ at the origin being a rational homology
  sphere), as well as an arithmetic one on the Milnor algebra of $f$,
  \cite{HollandMond-Scand98}. In \cite{Schulze-AdvGeom10}, using the
  Gau\ss--Manin connection, this is
  extended to a list of conditions on an isolated hypersurface
  singularity, each one of which forces the implication [$D$ has LCT]
  $\Rightarrow$ [$D$ is quasi-homogeneous].
\item For a version regarding more general connections, see
\cite{CalderonNarvaez-LMS09}.
\end{asparaenum}
\end{rmk}

A plane curve satisfies LCT if and only it is locally
quasi-homogeneous, \cite{CalderonCastroMondNarvaez-Helvetici02}.  By
\cite{CastroMondNarvaez-TAMS96}, free locally quasi-homogeneous
divisors satisfy LCT in any dimension. 
By
\cite{GrangerSchulze-Compos06}, in dimension three, free
divisors with LCT must be locally Euler-homogeneous.  
Conjecturally, LCT implies local Euler-homogeneity
\cite{CalderonCastroMondNarvaez-Helvetici02}. The converse is false,
see for example \cite{CastroUcha-PAMS05}.  
The classical example of
rotating lines with varying cross-ratio $f=xy(x+y)(x+yz)$ is free,
satisfies LCT and is locally Euler-homogeneous, but only weakly
quasi-homogeneous,
\cite{CalderonCastroMondNarvaez-Helvetici02}. 
In
\cite{CastroGagoHartilloUcha-Revista07}, the effect of the Spencer
property on LCT is discussed in the presence of homogeneity
conditions. 
For locally
quasi-homogeneous divisors (or if the non-free locus is
zero-dimensional), LCT implies $(B_1)$,
\cite{CastroUcha-Steklov02,Torrelli-survey}. In particular, LCT
implies $(B_1)$ for divisors with isolated singularities.  In
\cite{GrangerSchulze-Manuscripta06} quasi-homogeneity of isolated
singularities is characterized in terms of a map of local cohomology
modules of logarithmic differentials.

A free divisor is \emph{linear free} if the (free) module $\Der_X(-\log
f)$ has a basis of linear vector fields.  In
\cite{GrangerMondNietoSchulze-Fourier09}, linear free divisors in
dimension at most $4$ are classified, and for these divisors LCT holds
at least 
on global sections. In the process, it is shown that LCT is implied if
the Lie algebra of linear logarithmic vector fields is reductive. The
example of $n\times n$ invertible upper triangular matrices acting on
symmetric matrices \cite[Ex.~5.1]{GrangerMondNietoSchulze-Fourier09}
shows that LCT may hold without the reductivity assumption. Linear
free divisors appear naturally, for example in quiver representations
and in the theory of prehomogeneous vector spaces and castling
transformations
\cite{BuchweitzMond-LMS324,SatoKimura-Nagoya77,GrangerMondSchulze-LMS11}.
Linear freeness is related to unfoldings and Frobenius structures
\cite{GregorioMondSevenheck-Compos09}.

Denote by $\Der_{X,0}(-\log f)$ the derivations $\theta$ with
$\theta\bullet f=0$. In the presence of a global Euler-homogeneity $E$
on $Y$ there is a splitting $\Der_X(-\log f)\cong R\cdot E \oplus
\Der_{X,0}(-\log f)$. 
Reading derivations as
operators of order one,
$\Der_{X,0}(-\log f)\subseteq  \ann_D(f^s)$. We write $S$ for
$\gr_{(0,1)}(D)$; if $y_i$ is the symbol of $\del_i$ then we have
$S=R[y]$.
\begin{dfn}
The inclusion $\Der_{X,0}(-\log f)\into \ann_D(f^s)$, via the order filtration,
defines a subideal of $\gr_{(0,1)}(\ann_{D}(f^s))\subseteq
\gr_{(0,1)}(D)=S$ called the \emph{logarithmic ideal} $L_f$ of
$\Var(f)$.
\end{dfn}
Note that the symbols of $\Der_X(-\log f)$ are in the ideal $R\cdot y$
of height $n$.
\begin{dfn}\label{dfn-Koszul-free}
If $\Der_X(-\log f)$ has a generating set (as an
$R$-module) whose symbols form a regular sequence on $S$, then $Y$ is
called \emph{Koszul free}.
\end{dfn}
As $\Der_X(-\log f)$ has rank $n$, a
Koszul free divisor is indeed free.  Divisors in the plane
\cite{KSaito-Tokyo} and locally quasi-homogeneous free divisors
\cite{CalderonNarvaez-Compos02, CalderonNarvaez-Steklov02} are Koszul
free. In the case of normal crossings, this has been used to make
resolutions for $D[s]\bullet f^s$ and $D[s]/D[s](\ann_{D[s]} f^s,f)$,
\cite{GrosNarvaez-Padova00}. A way to distill invariants from
resolutions of $D[s]\bullet f^s$ is given in \cite{Arcadias-JPAA10}.
The logarithmic module $\tilde M^{\log f}=D/D\cdot \Der_X(-\log f)$ has
in the Spencer case (see \cite{CastroUcha-Steklov02,
  CalderonNarvaez-Fourier05}) a natural free resolution of Koszul
type.

For Koszul-free divisors, the ideal $D\cdot \Der_X(-\log f)$ is
holonomic \cite{Calderon-Ecole99}.  By
\cite[Thm.~7.4]{GrangerMondNietoSchulze-Fourier09}, in the presence of
freeness, the Koszul property is equivalent to the local finiteness of
Saito's logarithmic stratification. This yields an algorithmic way to
certify (some) free divisors as not locally quasi-homogeneous, since free
locally quasi-homogeneous divisors are Koszul free. Based on similar
ideas, one may devise a test for strong local Euler-homogeneity
\cite[Lem.~7.5]{GrangerMondNietoSchulze-Fourier09}.  See
\cite{Calderon-Ecole99} and \cite[\S2]{Torrelli-survey} for relations
of Koszul freeness to perversity of the logarithmic de Rham complex.

Castro-Jim\'enez and Ucha established conditions for $Y=\Var(f)$ to
have LCT in terms of $D$-modules
\cite{CastroUcha-JSC01,CastroUcha-Steklov02, CastroUcha-Exp04} for
certain free $f$.  For example, LCT is equivalent to $(A_1)$ for
Spencer free divisors. Calder\'on-Moreno and Narv\'aez-Macarro
\cite{CalderonNarvaez-Fourier05} proved that free divisors have LCT if
and only if the natural morphism
$\calD_X\otimes_{V^0(\calD_X)}^L\calO_X(Y)\to \calO_X(*Y)$ is a
quasi-isomorphism, $\calO_X(Y)$ being the meromorphic functions with simple
pole along $f$.  
For Koszul free $Y$, one has
at least $\calD_X\otimes_{V^0(\calD_X)}^L\calO_X(Y)\cong
\calD_X\otimes_{V^0(\calD_X)}\calO_X(Y)$. A similar condition ensures
that the logarithmic de Rham complex is perverse
\cite{Calderon-Ecole99,CalderonNarvaez-Fourier05}. The two results are
related by duality between logarithmic connections on $\calD_X$ and
the $V$-filtration
\cite{CastroUcha-Steklov02,CalderonNarvaez-Fourier05,CastroUcha-Paris04}.

It is unknown how LCT is related to $(A_1)$ in general, but for
quasi-homogeneous polynomials with isolated singularities the two
conditions are equivalent,
\cite{Torrelli-survey}.

\subsection{Logarithmic linearity}
\begin{dfn}
We say that $f\in R$ satisfies $(L_s)$ if the characteristic ideal of
$\ann_D(f^s)$ is generated by symbols of derivations.
\end{dfn}
Condition $(L_s)$  holds for isolated
singularities
\cite{Yano}, locally
quasi-homogeneous free divisors \cite{CalderonNarvaez-Compos02}, and
locally strongly Euler-homogeneous holonomic tame divisors \cite{Walther-zeta}.
Also, $(L_s)$ plus $(B_1)$ yields $(A_1)$ for locally Euler-homogeneous $f$ by
\cite{Kashiwara-bfu}, see \cite{Torrelli-survey}.

The logarithmic ideal supplies an interesting link between
$\Omega^\bullet_X(\log f)$ and $\ann_D(f^s)$ via approximation complexes:
if $f$ is holonomic, strongly locally Euler-homogeneous and also tame
then the complex $(\Omega^\bullet_X(\log f)[y],\ydx)$ is a resolution of
the logarithmic ideal $L_f$, and $S/L_f$ is a Cohen--Macaulay domain of
dimension $n+1$; if $f$ is in fact free, $S/L_f$ is a complete
intersection \cite{Narvaez-Contemp08,Walther-zeta}.
\begin{que}
For locally Euler-homogeneous divisors, is $\ann_D(f^s)$ related to
the cohomology of $(\Omega^\bullet_X(\log f)[y],\ydx)$?
\end{que}

\section{Characteristic variety}

We continue to assume that $X=\CC^n$.
For $f\in R$ let $U_f$ be the open set defined by  $\de f\not=0\not
=f$. 
Because of the functional equation,
$\calM_f(s)$ is coherent over $D$ \cite{Bernstein-Nauk72,Kashiwara-bfu},
and the restriction of 
$\charVar(D[s]\bullet f^s)$ to $U_f$ is
\begin{equation}\label{eqn-charvar-U}
\overline{\left\{(\xi,s\frac{\de f(\xi)}{f(\xi)})\mid
  \xi \in U_f, s\in\CC\right\}}^{\text{Zariski}},
\end{equation}
an $(n+1)$-dimensional
involutive subvariety of $T^*{U_f}$,
\cite{Kashiwara-AMSbook}. Ginsburg \cite{Ginsburg-Invent86} gives a formula
for the characteristic cycle of $D[s]\bullet mf^s$ in terms of an
intersection process for holonomic sections $m$.

In favorable cases, more can be said.
By \cite{CalderonNarvaez-Compos02}, if the
divisor is reduced, free and locally quasi-homogeneous then
$\ann_{D[s]}(f^s)$ is generated by derivations, both $\calM_f(s)$ and
$\calN_f(s)$ have Koszul--Spencer type resolutions, and in particular the
characteristic varieties are complete intersections. In the more general
case where
$f$ is holonomic, locally strongly Euler-homogeneous and tame,
$\ann_D(f^s)$ is still generated by
order one operators and the ideal of symbols of $\ann_D(f^s)$ (and hence the
characteristic ideal of $\calM_f(s)$ as well)
is a Cohen--Macaulay prime ideal, \cite{Walther-zeta}. Under these
hypotheses, the characteristic ideal of $\calN_f(s)$ is Cohen--Macaulay
but usually not prime.

\subsection{Stratifications}
By \cite{KashiwaraSchapira-Acta79}, the resolution theorem
of Hironaka can be used to show that there is a stratification of
$\CC^n$ such that for each holonomic $D$-module $M$,
$\charCycle(M)=\bigsqcup_{\sigma\in\Sigma}
\mu(M,\sigma)T^*_\sigma$ where $T^*_\sigma$ is the closure of the conormal bundle
of the smooth stratum $\sigma$ in $\CC^n$ and $\mu(M,\sigma)\in\NN$.

For $D[s]\bullet f^s/D[s]\bullet f^{s+1}$ Kashiwara proved that if one
considers a Whitney stratification $S$ for $f$ (for example the
``canonical'' stratification in \cite{DamonMond-Invent91}) then the
characteristic variety of the $D$-module $\calN_f(s)$ is in the union
of the conormal varieties of the strata $\sigma\in S$, \cite{Yano}.

If one slices a pair $(X,D)$ of a smooth
space and a divisor with a hyperplane, various invariants of the
divisor will behave well provided that the hyperplane is not
``special''. A prime example are Bertini and Lefschetz theorems. For
$D$-modules, Kashiwara defined the notion of \emph{non-characteristic
  restriction}: the smooth hypersurface $H$ is non-characteristic for
the $D$-module $M$ if it meets each component of the characteristic
variety of $M$ transversally (see  \cite{Pham-book79} for an
exposition). The condition assures that the inverse image functor
attached to the embedding $H\into X$ has no higher derived functors
for $M$. In \cite{DimcaMaisonobeSaitoTorrelli-MathAnn06} these ideas
are used to show that the $V$-filtration, and hence the multiplier
ideals as well as nearby and vanishing cycle sheaves, behave nicely
under non-characteristic restriction.

\subsection{Deformations}

Varchenko proved, via establishing constancy of Hodge numbers, that in a
$\mu$-constant family of isolated singularities, the spectrum is
constant \cite{Varchenko-Nauk82}.
In \cite{DimcaMaisonobeSaitoTorrelli-MathAnn06} it is shown that the
formation of the spectrum along the divisor $Y\subseteq X$ commutes with the
intersection with a hyperplane transversal to any stratum of a Whitney
regular stratification of $D$,
and a weak
generalization of Varchenko's constancy results for certain
deformations of non-isolated singularities is derived.

In contrast, the Bernstein--Sato polynomial may not be constant along
$\mu$-constant deformations.  Suppose $f(x)+\lambda g(x)$ is a
$1$-parameter family of plane curves with isolated singularities at
the origin. If the Milnor number $\dim_\CC(R/J_{(f+\lambda g)})$ is
constant in the family, the singularity germs in the family are
topologically equivalent \cite{LeRamanujam-AJM76}; for discussion see
\cite[\S2]{Dimca-SingTopHyp}. However, in such a family $b_f(s)$ can
vary, as it is a differential invariant. Indeed, $f+\lambda 
g=x^4+y^5+\lambda xy^4$ has constant Milnor number 20, but the general
curve (not quasi-homogeneous in any coordinate system, as
$\rho_{f+\lambda g}$ is not symmetric about $-1$, see Subsection
\ref{subsec-bfu} above) has $-\rho_{f+\lambda g}=\{1\}\cup
\frac{1}{20}\{9,11,13,14,17,18,19,21,22,23,26,27\}$ while the special
curve has $-\rho_{f}=-\rho_{f+\lambda g}\cup\{-31/20\}\smallsetminus
\{-11/20\}$. See \cite{CassouNogues-Fourier86} for details and similar
examples based on Newton polytope considerations, and
\cite{Stahlke-thesis} for deformations of plane diagonal curves.

\section{Milnor fiber and monodromy}

\subsection{Milnor fibers}\label{sec-Milnor}
Let $B(p,\eps)$ denote the $\eps$-ball around $p\in\Var(f)\subseteq
\CC^n$. Milnor \cite{Milnor} proved that the diffeomorphism type of
the open real manifold
\[
M_{p,t_0,\eps}=B(p,\eps)\cap \Var(f-t_0)
\]
is independent of $\eps,t_0$ as long as $0<|t_0|\ll\eps\ll 1$. For
$0<\tau\ll\eps\ll 1$ denote by
$M_p$ the fiber of the bundle $B(p,\eps)\cap \{q\in\CC^n\mid
0<|f(q)|<\tau\}\to f(q)$.

The direct image functor for $D$-modules to the projection
$\CC^n\times \CC\to \CC$, $(x,t)\mapsto t$ turns
the $D_{x,t}$-module $\calB_f$ into the
\emph{Gau\ss\--Manin system} $\calH_f$.
The $D$-module restriction of $H^k(\calH_f)$ to
$t=t_0$ is the $k$-th cohomology of the Milnor fibers along
$\Var(f)$ for $0<|t_0|<\tau$.

Fix a $k$-cycle $\sigma\in H_p(\Var(f-t_0))$ and choose $\eta\in
H^k(\calH_f)$. Deforming $\sigma$ to a $k$-cycle over $t$ using the
Milnor fibration, one can evaluate $\int_{\sigma_t}\eta$.  The
Gau\ss--Manin system has Fuchsian singularities and these periods are
in the Nilsson class \cite{Malgrange-integr74}.
For example,  the classical Gau\ss\ hypergeometric
function saw the light of day the first time as solution to a system
of differential equations attached to the variation of the Hodge
structure on an elliptic curve (expressed as integrals of the first
and second kind) \cite{BrieskornKnoerrer}. In \cite{Pham-book79} this
point of view is taken to be the starting point. The techniques
explained there form the foundation for many connections between
$f^s$ and singularity invariants attached to $\Var(f)$.

In \cite{Budur-MathAnn03}, a bijection (for $0<\alpha\le 1$) is
established between a subset of the jumping numbers of $f$ at $p\in\Var(f)$ and
the support of the \emph{Hodge spectrum} \cite{Steenbrink-Asterisque89}
\[
\Sp(f)=\sum_{\alpha\in\QQ}n_\alpha(f)t^\alpha,
\]
with $n_\alpha(f)$ determined by the size of
the $\alpha$-piece of Hodge component of the cohomology of the Milnor
fiber of $f$ at $p$. See also
\cite{Saito-MathAnn93,Varchenko-Izvestiya81}, and
\cite{Steenbrink-survey} for a survey on Hodge invariants.
We refer to \cite{Budur-survey,Saito-Tokyo09} for many more
aspects of this part of the story.

\subsection{Monodromy}
The vector spaces
  $H^k(M_{p,t_0,\eps},\CC)$ form a smooth vector bundle over a
  punctured disk $\CC^*$. The linear transformation $\mu_{f,p,k}$
  on $H^k(M_{p,t_0,\eps},\CC)$ induced by
  $p\mapsto p\cdot \exp(2\pi i \lambda)$ is the $k$-th monodromy of $f$ at
  $p$. Let $\chi_{f,p,k}(t)$ denote the characteristic polynomial of
  $\mu_{f,p,k}$, set
\[
e_{f,p,k}=\{\gamma\in\CC\mid \gamma\text{ is an eigenvalue of
}\mu_{f,p,k}\}
\]
and put $e_{f,p}=\bigcup e_{f,p,k}$.

For most (in a quantifiable sense) divisors $f$ with given Newton
diagram, a combinatorial recipe can be given that determines the
alternating product $\prod (\chi_{f,p,k}(t))^{(-1)^k}$
\cite{Varchenko-Invent76}, similarly to A'Campo's formula in terms of
an embedded resolution \cite{ACampo-Helvetici75}.

\subsection{Degrees, eigenvalues, and Bernstein--Sato polynomial}

By \cite{Malgrange-evan, Kashiwara-LNM1016}, the exponential function maps
the root set of the local analytic Bernstein--Sato polynomial of $f$
at $p$ onto $e_{f,p}$. The set $\exp(-2\pi i\tilde\rho_{f,p})$ is the set
of eigenvalues of the monodromy on the Grothendieck--Deligne
vanishing cycle sheaf $\phi_f(\CC_{X,p})$. This was shown in
\cite{Saito-Bull94} by algebraic microlocalization.

If $f$ is an isolated singularity, the Milnor fiber $M_f$ is a bouquet
of spheres, and $H^{n-1}(M_f,\CC)$ can be identified with the Jacobian
ring $R/J_f$ as vector space. Moreover, if $f$ is quasi-homogeneous,
then under this identification $R/J_f$ is a $\QQ[s]$-module, $s$
acting via the Euler operator, and $\tilde\rho_f$ is in bijection with
the degree set of the nonzero quasi-homogeneous elements in
$R/J_f$. For non-isolated singularities, most of this breaks down,
since $R/J_f$ is not Artinian in that case.  However, for homogeneous
$f$, consider the \emph{Jacobian module}
\[
H^0_\frakm(R/J_f)=\{g+J_f\mid \exists k\in\NN, \forall i, x_i^kg\in J_f\}
\]
and
the canonical $(n-1)$-form
\[
\eta=\sum_ix_i\de x_1\wedge\cdots\wedge\widehat{\de
  x_i}\wedge\cdots\wedge \de x_n.
\]
Every class in $H^{n-1}(M_f;\CC)$ is of the form $g\eta$ for suitable
$g\in R$, and there is a filtration on $H^{n-1}(M_f,\CC)$ induced by
integration of $\calB_f$ along $\del_1,\ldots,\del_n$, with the
following property: if $g\in R$ is the smallest degree homogeneous
polynomial such that $g\eta$ represents a chosen element of
$H^{n-1}(M_f,\CC)$ then $b_f(-(\deg(g\eta))/\deg(f))=0$,
\cite{Walther-Bernstein}.  Suppose the projective variety defined by
$f$ has isolated singularities. Then by
\cite{Saito-0602527,Walther-zeta}, with $1\le k\le d$ and
$\lambda=\exp(2\pi \sqrt{-1}k/d)$, the following holds: $\dim_\CC
     [H^0_\frakm(R_n/\Jac(f))]_{d-n+k}\le \dim_\CC \gr^{{\rm
         Hodge}}_{n-2}(H^{n-1}(M_f,\CC)_\lambda)$ where the right hand
     side indicates the $\lambda$-eigenspace of the associated graded
     object to the Hodge filtration on $H^{n-1}(M_{f})$.

\subsection{Zeta functions}\label{subsec-zeta}
The zeta function $Z_f(s)$ attached to a divisor $f\in R$ is the rational
function
\[
Z_f(s)=\sum_{I\subseteq S}\chi(E^*_I)\prod_{i\in
  I}\frac{1}{N_is+\nu_i}
\]
where $\pi\colon (Y,\bigcup_I E_i)\to (\CC^n,\Var(f))$ is an embedded
resolution of singularities, and $N_i$ (resp.\ $\nu_i-1$) are the
multiplicities of $E_i$ in $\pi^*(f)$ (resp.\ in the Jacobian of
$\pi$). By results of Denef and Loeser \cite{DenefLoeser-JAMS92},
$Z_f(s)$ is independent of the resolution.

\begin{cnj}[Topological Monodromy Conjecture]
\begin{itemize}~
\item[(SMC)] Any pole of $Z_f(s)$ is a root of the
Bernstein--Sato polynomial $b_f(s)$.
\item[(MC)] Any pole of $Z_f(s)$ yields under exponentiation an
  eigenvalue of the monodromy operator at some $p\in\Var(f)$.
\end{itemize}
\end{cnj}
The strong version (SMC) implies (MC) by
\cite{Malgrange-isolee,Kashiwara-LNM1016}.
Each version allows a generalization to ideals. 

(SMC), formulated by Igusa \cite{Igusa-AMS14} and Denef--Loeser
\cite{DenefLoeser-JAMS92} holds for
\begin{itemize}
\item reduced curves by \cite{Loeser-AJM88}  with a
  discussion on the nature of the poles 
  by Veys \cite{Veys-Bull93,Veys-JLMS90,Veys-Manuscripta95};
\item certain Newton-nondegenerate divisors by  \cite{Loeser-Crelle90};
\item some hyperplane arrangements (see Section \ref{sec-arr});
\item  monomial ideals in any dimension
   by \cite{HowaldMustataYuen-PAMS07}.
\end{itemize}
Additionally, Conjecture (MC) holds for
\begin{itemize}
\item  bivariate ideals by Van Proeyen and Veys
  \cite{vanProeyenVeys-Fourier10};
\item all hyperplane arrangements by  \cite{BudurMustataTeitler-Geo11,BudurSaitoYuzvinsky-JLMS11};
\item some partial cases:
  \cite{ArtalCassouLuengoMelle-Ecole,LemahieuVeys-IMRN} some surfaces;
  \cite{ArtalCassouLuengoMelle-Memoirs} quasi-ordinary power series;
  \cite{LichtinMeuser-Compos85,Loeser-Crelle90} in certain Newton
  non-degenerate cases; \cite{Igusa-AJM92,KimuraSatoZhu-AJM90}
  for invariants of prehomogeneous vector spaces;
  \cite{LemahieuVanProeyen-TAMS} for  nondegenerate surfaces.
\end{itemize}
Strong evidence for (MC) for $n=3$ is procured in
\cite{Veys-Crelle06}.  The articles
\cite{Rodrigues-Cambridge04,NemethiVeys-GeomTop12} explore what (MC)
could mean on a normal surface as ambient space and gives some results
and counterexamples to naive generalizations.  See also
\cite{Denef-Bourbaki} and the introductions of
\cite{Bories-thesis,BoriesVeys-1306.6012} for more details in survey
format.

\medskip

A root of $b_f(s)$, a monodromy eigenvalue, and a pole of $Z_f(s)$ may
have multiplicity; can the monodromy conjecture be strengthened to
include multiplicities?  This version of (SMC) was proved for reduced
bivariate $f$ in \cite{Loeser-AJM88}; in
\cite{MelleTorrelliVeys-JT09,MelleTorrelliVeys-JA10} it is proved for
certain nonreduced bivariate $f$, and for some trivariate ones.

A different variation, due to Veys, of the conjecture is the
following. Vary the definition of $Z_f(s)$ to $Z_{f;g}(s)=\sum_{I\subseteq
  S}\chi(E^*_I)\prod_{i\in I}\frac{1}{N_is+\nu'_i}$ where $\nu'_i$ is
the multiplicity of $E_i$ in the pullback along $\pi$ of some
differential form $g$. (The standard case is when $g$ is the volume
form). Two questions arise: (1) varying over a suitable set $G$ of forms $g$,
can one generate all roots of $b_f(s)$ as poles of the resulting zeta
functions?
And if so, can one (2) do this such that the pole sets of all zeta
functions so constructed are always inside $\rho_f$, so that
\[
\rho_f=\{\alpha\mid \exists g\in G,
\lim_{s\to\alpha}Z_{f;g}(s)=\infty\}\quad ?
\]
N\'emethi and Veys \cite{NemethiVeys-BLMS10,NemethiVeys-GeomTop12}
prove a weak version: if $n=2$ then monodromy eigenvalues are
exponentials of poles of zeta functions from differential forms.

The following is discussed in \cite{Bories-zeta-bfu}. For some ideals
with $n=2$, (1) is false for the topological
zeta function (even for divisors: consider
$xy^5+x^3y^2+x^4y$). For monomial ideals with
two generators in $n=2$, (1) is correct; with more than two generators
it can fail. Even in the former  case, (2) can be
false.

\section{Multi-variate versions}

If $f=(f_1,\ldots,f_r)$ defines a map $f\colon \CC^n\to\CC^r$, several
$b$-functions can be defined:
\begin{asparaenum}
\item The univariate Bernstein--Sato polynomial $b_f(s)$ attached to the ideal
  $(f)\subseteq R$ from \cite{BudurMustataSaito-Compos06}.
\item The multi-variate Bernstein--Sato polynomials $b_{f,i}(s)$
  of all $b(s)\in\CC[s_1,\ldots,s_r]$ such that there is an
  equation $P(x,\del,s)\bullet f_if^s=b(s)f^s$ in multi-index notation.
\item The multi-variate Bernstein--Sato ideal $B_{f,\mu}(s)$ for
  $\mu\in\NN^r$
  of all $b(s)\in\CC[s_1,\ldots,s_r]$ such that there is an
  equation $P(x,\del,s)\bullet f^{s+\mu}=b(s)f^s$ in multi-index
  notation. The most interesting case is $\mu=\boldone=(1,\ldots,1)$.
\item  The multi-variate Bernstein--Sato ideal $B_{f,\Sigma}(s)$ 
  of all $b(s)\in\CC[s_1,\ldots,s_r]$ that
  multiply $f^s$ into  $\sum D[s]f_i f^{s}$ in multi-index notation.
\end{asparaenum}

The Bernstein--Sato polynomial in (1) above has been studied in the
case of a monomial ideal in \cite{BudurMustataSaito-CIA06} and more
generally from the point of view of the Newton polygon in
\cite{BudurMustataSaito-MRL06}. While the roots for monomial ideals do
not depend just on the Newton polygon, their residue classes modulo
$\ZZ$ do.

\medskip

Non-triviality of the quantities in (2)-(4) have been established in
\cite{Sabbah-Compos87,Sabbah-Compos87II,Sabbah-ecole87}, but see also
\cite{Bahloul-Compos05}.
The ideals $B_{f,\mu}(s)$ and $B_{f,\Sigma}(s)$ do not have to be
principal, \cite{CastroUcha-JSC04,BahloulOaku-JSC10}.
In \cite{Sabbah-Compos87,Gyoja-Kyoto93} it is
shown that $B_{f,\mu}(s)$ contains a polynomial that factors into linear
forms with non-negative rational coefficients and positive constant term.
Bahloul and Oaku \cite{BahloulOaku-JSC10} show that these ideals are
local in the sense of \eqref{eq-local-bfu}.

The following would generalize Kashiwara's result in the univariate
case as well as the results of Sabbah and Gyoja above.
\begin{cnj}[\cite{Budur-multi12}]\label{cnj-Sabbah}
The Bernstein--Sato ideal $B_{f,\mu}(s)$ is
is generated by products of linear
forms $\sum \alpha_is_i+a$ with $\alpha_i,a$ non-negative rational and
$a>0$.
\end{cnj}
For $n=2$, partial results by Cassou-Nogu\`es and Libgober exist
\cite{CassouNoguesLibgober-Knot11}. In \cite{Budur-multi12} it is
further conjectured that the Malgrange--Kashiwara result,
exponentiating $\rho_{f,p}$ gives $e_{f,p}$, generalizes: monodromy
in this case is defined in \cite{Verdier-Asterisque83}, and Sabbah's
specialization functor $\psi_f$ from \cite{Sabbah-Duke90} takes on the
r\^ole of the nearby cycle functor, and conjecturally exponentiating
the variety of $B_{f,p}(s)$ yields the uniform support (near $p$) of Sabbah's
functor. The latter conjecture would imply Conjecture
\ref{cnj-Sabbah}.

Similarly to the one-variable case, if $V(n,d,m)$ is the vector space
of (ordered) $m$-tuples of polynomials in $x_1,\ldots,x_n$ of degree
at most $d$, there is an algebraic stratification of $V(n,d,m)$ such
that on each stratum the function $V\ni f=(f_1,\ldots,f_m)\mapsto
b_f(s)$ is constant. Corresponding results for the Bernstein--Sato
ideal $B_{f,\boldone}(s)$ hold by
\cite{BrianconMaisonobe-Tokyo00}.

\section{Hyperplane arrangements}\label{sec-arr}

A \emph{hyperplane arrangement} is a divisor of the form
\[
\calA=\prod_{i\in I}\alpha_i
\]
where each $\alpha_i$ is a polynomial of degree one. We denote
$H_i=\Var(\alpha_i)$. Essentially all
information we are interested in is of local nature, so we assume that
each $\alpha_i$ is a form so that $\calA$ is
\emph{central}. If there is a coordinate change in $\CC^n$ such that
$\calA$ becomes the product of polynomials in disjoint sets of variables,
the arrangement is \emph{decomposable}, otherwise it is
\emph{indecomposable}.

A \emph{flat} is any (set-theoretic) intersection $\bigcap_{i\in J}H_i$
where $J\subseteq I$.
The \emph{intersection lattice} $L(\calA)$ is the partially ordered set
consisting of the collection of
all flats, with order given by inclusion.

\subsection{Numbers and parameters}
Hyperplane arrangements satisfy $(B_1)$ everywhere
\cite{Walther-Bernstein}.
Arrangements satisfy $(A_1)$ everywhere if they
decompose into a union of a generic and a hyperbolic arrangement
\cite{Torrelli-Bull04}, and if they are tame
\cite{Walther-zeta}.  Terao conjectured that all hyperplane
arrangements satisfy $(A_1)$; some of them fail $(A_s)$, \cite{Walther-zeta}.

\medskip

Apart from recasting various of the previously encountered problems in
the world of arrangements, a popular study is the following: choose a
discrete invariant $I$ of a divisor. Does the function $\calA\mapsto
I(\calA)$ factor through the map $\calA\mapsto L(\calA)$?  Randell
showed that if two arrangements are connected by a one-parameter
family of arrangements which have the same intersection 
lattice, the complements are
diffeomorphic \cite{Randell-PAMS89} and the isomorphism type of
the Milnor fibration is constant
\cite{Randell-PAMS97}. Rybnikov
\cite{Rybnikov-Nauk11,ArtalCarmonaCogolludoMarco-PureMath06} showed on
the other hand that there are arrangements (even in the projective plane)
with equal lattice but different complement. In particular, not all
isotopic arrangements can be linked by a smooth deformation.

\subsection{LCT and logarithmic ideal}
The most prominent positive result is by Brieskorn: the
\emph{Orlik--Solomon algebra} $\OS(\calA)\subseteq
\Omega^\bullet(\log\calA)$ generated by the forms $\de \alpha_i/\alpha_i$ is
quasi-isomorphic to $\Omega^\bullet(*\calA)$, hence to the singular
cohomology algebra of $U_\calA$, \cite{Brieskorn-tresses}. The
relation with combinatorics was given in 
\cite{OrlikSolomon-Invent80,OrlikTerao-book}. For a survey on the
Orlik--Solomon algebra, see \cite{Yuzvinsky-OS01}. The best known open
problem in this area is
\begin{cnj}[\cite{Terao-Kyoto78}]
$\OS(\calA)\to
\Omega^\bullet(\log\calA)$ is a quasi-isomorphism.
\end{cnj}
While the general case remains open, Wiens and Yuzvinsky
\cite{WiensYuzvinsky-TAMS97} proved it for tame arrangements, and also
if $n\le 4$. The techniques are based on
\cite{CastroMondNarvaez-TAMS96}.

\subsection{Milnor fibers}
There is a survey article by Suciu on complements, Milnor fibers, and
cohomology jump loci \cite{Suciu-1301.4851}, and \cite{Budur-survey}
contains further information on the topic. It is not known whether
$L(\calA)$ determines the Betti numbers (even less the Hodge numbers) of
the Milnor fiber of an arrangement. The first Betti number of the
Milnor fiber $M_\calA$ at the origin is stable under intersection with a generic
hyperplane (if $n>2$). But it is unknown whether the first Betti
number of an arrangement in $3$-space is a function of the lattice
alone. By \cite{Dimcaetal-1305.5092}, this is so for collections of up
to 14 lines with up to $5$-fold intersections in the projective
plane. See also \cite{Libgober-Sapporo12} for the origins of the
approach. By \cite{BudurDimcaSaito-AMS11}, a lower combinatorial bound
for the $\lambda$-eigenspace of $H^1(M_\calA)$ is given under
favorable conditions on $L$. If $L$ satisfies stronger conditions, the
bound is shown to be exact. In any case,  \cite{BudurDimcaSaito-AMS11}
gives an algebraic, although
perhaps non-combinatorial, formula for the Hodge pieces in terms of
multiplier ideals.

By \cite{OrlikRandell-Arkiv93}, the Betti numbers of $M_\calA$ are
combinatorial if $\calA$ is generic. See also \cite{CohenSuciu-JLMS95}.

\subsection{Multiplier ideals}
\mustata\ gave a formula for the multiplier ideals of arrangements,
and used it to show that the log-canonical threshold is a function of
$L(\calA)$. The formula is somewhat hard to use for showing that each
jumping number is a lattice invariant; this problem was solved in
\cite{BudurSaito-MathAnn10}.  Explicit formulas in low
dimensional cases follow from the spectrum formulas given
there and in \cite{Yoon-1211.1689}. Teitler \cite{Teitler-PAMS08} improved
\mustata's formula for multiplier ideals to not necessarily reduced
hyperplane arrangements \cite{Mustata-TAMS06}.

\subsection{Bernstein--Sato polynomials}
By \cite{Walther-Bernstein},
$\rho_\calA\cap\ZZ=\{-1\}$; by \cite{Saito-0602527}, $\rho_\calA\subseteq(-2,0)$.
There are few  classes of arrangements with explicit formul\ae\ for
their  Bernstein--Sato
polynomial:
\begin{itemize}
\item  Boolean (a normal crossing arrangement, locally given by
$x_1\cdots x_k$);
\item hyperbolic (essentially an arrangement in two
variables);
\item generic (central, and all intersections of $n$ hyperplanes equal
  the origin).
\end{itemize}
The first case is trivial, the second is easy, the last is
\cite{Walther-Bernstein} with assistance from \cite{Saito-Compos07}.
Some interesting computations are in
\cite{BudurSaitoYuzvinsky-JLMS11}, and
\cite{Budur-multi12} has a partial confirmation of the multi-variable
Kashiwara--Malgrange theorem.
The Bernstein--Sato polynomial is not determined by the intersection
lattice, \cite{Walther-zeta}.

\subsection{Zeta functions}
Budur, \mustata\ and  Teitler \cite{BudurMustataTeitler-Geo11} show:
(MC) holds for arrangements, and
in order to prove (SMC), it suffices to show
the following conjecture.
\begin{cnj}\label{n/d}
For all indecomposable central arrangements with $d$
planes in $n$-space,
$b_\calA(-n/d)=0$.
\end{cnj}
The idea is to use the resolution of singularities 
obtained by blowing up the dense
edges from \cite{SchechtmanTeraoVarchenko}. The corresponding
computation of the zeta function is inspired from
\cite{Igusa-Crelle74,Igusa-Crelle75}. The number $-n/d$ does not have to be
the log-canonical threshold.
By \cite{BudurMustataTeitler-Geo11},
Conjecture~\ref{n/d} holds in a number of cases, including reduced
arrangements in dimension $3$. By \cite{Walther-zeta} it holds for
tame arrangements.

Examples of Veys (in $4$ variables) show that (SMC) may hold even if
Conjecture~\ref{n/d} were false in general, since $-n/d$ is not always
a pole of the zeta function
\cite{BudurSaitoYuzvinsky-JLMS11}. However, in these examples, $-n/d$
is in fact a root of $b_f(s)$.

For arrangements, each monodromy eigenvalue can be captured by zeta
functions in the sense of N\'emethi and Veys, see Subsection
\ref{subsec-zeta}, but not necessarily all of $\rho_\calA$ (Veys and
Walther, unpublished).

\section{Positive characteristic}
Let here $\FF$ denote a field of characteristic $p>0$. The theory of
$D$-modules is rather different in positive characteristic compared to
their behavior over the complex numbers. There are several reasons for
this:
\begin{enumerate}
\item On the downside, the ring $D_p$ of $\FF$-linear differential operators
  on $R_p=\FF[x_1,\ldots,x_n]$ is no longer finitely generated: as an
  $\FF$-algebra it is generated by the elements $\del^{(\alpha)}$,
  $\alpha\in\NN^n$, which act via
  $\del^{(\alpha)}\bullet(x^\beta)={\beta\choose\alpha}x^{\beta-\alpha}$.
\item As a trade-off, one has access to the Frobenius morphism
  $x_i\mapsto x_i^p$, as well as the Frobenius
  functor $F(M)=R'\otimes_R M$ where $R'$ is the $R-R$-bimodule on
  which $R$ acts via the identity on the left, and via the Frobenius on the
  right. Lyubeznik \cite{Lyubeznik-Fmods} created the category of
  $F$-finite $F$-modules and proved striking finiteness results. The
  category includes many interesting $D_p$-modules, and
  all $F$-modules are $D_p$-modules.
\item Holonomicity is more complicated, see
  \cite{Bogvad-02}.
\end{enumerate}
A most surprising consequence of Lyubeznik's ideas is that in positive
characteristic the property $(B_1)$ is meaningless: it holds for every
$f\in R_p$, \cite{AlvarezBlickleLyubeznik-LRM05}. The proof uses in
significant ways the difference between $D_p$ and the Weyl algebra.
In particular, the theory of Bernstein--Sato polynomials is rather
different in positive characteristic. In \cite{Mustata-JA09}
a sequence of Bernstein--Sato polynomials is attached to a
polynomial $f$ assuming that the Frobenius morphism is finite on $R$
(\emph{e.g.}, if $\FF$ is finite or algebraically closed); these
polynomials are then linked to test ideals,
the finite characteristic counterparts to multiplier ideals. In
\cite{BlickleMustataSmith-TAMS09} variants of our modules
$\calM_f(\gamma)$ are introduced and 
\cite{NunezPerez-1302.3327} shows that simplicity of these modules
detects the $F$-thresholds from
\cite{MustataTakagiWatanabe-ECM05}. These are
cousins of the jumping numbers of multiplier ideals and 
related to the Bernstein--Sato polynomial via base-$p$-expansions; see also
\cite{WittHernandez-13}. The Kashiwara--Brylinski intersection homology
module was shown to exist in positive
characteristic by Blickle in his thesis, \cite{Blickle-MathAnn04}.

\section{Appendix: Computability (by A.~Leykin)}

Computations around $f^s$ can be carried out by hand in special
cases. Generally, the computations are enormous and computers are
required (although not often sufficient).  One of the earliest such
approaches are in \cite{BrianconGrangerMaisonobeMiniconi-Fourier89,
  AleksandrovKistlerov-Contemp92}, but at least implicitly
Buchberger's algorithm in a Weyl algebra was discussed as early as
\cite{Castro-thesis}. An algorithmic approach to the isolated
singularities case~\cite{Maisonobe-LMS94} preceded the general
algorithms based on \GBs\ in a non-commutative setting outlined below.

\subsection{\GBs}
The \emph{monomials} $x^\alpha\p^\beta$ with $\alpha,\beta\in\NN^n$
form a $\CC$-basis of $D$; expressing $p\in D$ as linear combination
of monomials leads to its {\em normal form}.  The monomial orders on
the commutative monoid $\monoid{x, \p}$ for which for all $i\in [n]$
the leading monomial of $\p_ix_i=x_i\p_i+1$ is $x_i\p_i$, can be used
to run Buchberger's algorithm in $D$. Modifications are needed in
improvements that exploit commutativity, but the na\"ive
Buchberger's algorithm works without any changes. See
\cite{KandriRodyWeispfenning-JSC90} for more general settings in
polynomial rings of solvable type.
Surprisingly, the
worst case complexity of \GBs\ computations in Weyl
algebras is {\em not} worse than in the commutative polynomial case:
it is doubly exponential in the number of
indeterminates~\cite{AschLeykin:GBB,Grigoriev-Chistov:D-complexity}.

\subsection{Characteristic variety}
A weight vector $(u,v)\in \ZZ^n\times \ZZ^n$ with
$u+v\geq 0$ induces a filtration of $D$,
\[
F_i = \Span{\CC}{x^\alpha \p^\beta \mid u\cdot\alpha + v\cdot\beta
  \leq i},\quad i\in \ZZ.
\]
The $(u,v)$-\emph{Gr\"obner deformation} of a left ideal  $I\subseteq D$
is
\[
\ini_{(u,v)}(I) = \Span{\CC}{\ini_{(u,v)}(P)\mid P\in I} \subseteq
\gr_{(u,v)}D,
\]
the ideal of {\em initial forms} of elements of $I$ with respect to
the given weight in the associated graded algebra.  It is possible to
compute \Grob\ deformations in the homogenized Weyl
algebra
\[
D^h = D\<h\>/\ideal{\p_ix_i-x_i\p_i-h^2,
    x_ih - hx_i, \p_ih - h\p_i, \mid
    1\le i\le n}
\]
see \cite{CastroNarvaez-preprint,OakuTakayama-JSC01}.
\Grob\ deformations
are the main topic of~\cite{HyperGeomBook}.

%
%

\subsection{Annihilator}

Recall the construction appearing in the beginning of
\S\ref{sec-Milnor}: $D_{x,t}$ acts on $D[s]f^s$; in
particular, the operator $-\p_t t$ acts as multiplication by
$s$. It is this approach that lead Oaku to an algorithm for
$\ann_{D[s]}(f^s)$, $\ann_D(f^s)$ and $b_f(s)$,
\cite{Oaku-local-b}. We outline the ideas.

Malgrange observed that
\begin{align}\label{eq:ann}
\ann_{D[s]}(f^s) &= \ann_{D_{x,t}}(f^s) \cap D[s],\\
\text{with}\quad \ann_{D_{x,t}}(f^s) &= \ideal{t-f,
  \p_1+\textstyle\frac{\p f}{\p x_1}\p_t, \ldots,
  \p_n+\textstyle\frac{\p f}{\p x_n}\p_t} \subseteq D_{x,t}. \label{eq:ann-Dxt}
\end{align}
The former can be found from the latter by eliminating $t$ and $\p_t$
from the ideal
\begin{equation}\label{eq:ideal-with-s-t-dt}
\ideal{s+t\p_t} + \ann_{D_{x,t}}(f^s) \subseteq D_{x,t}\<s\>;
\end{equation}
of course $s=-\del_tt$ does not commute with $t,\del_t$ here.

Oaku's method for $\ann_{D[s]}(f^s)$ accomplished the elimination
by augmenting two commuting indeterminates:
\begin{align}\label{eq:ann-Oaku}
\begin{split}
\ann_{D[s]}(f^s) &= I'_f \cap D[s],\\ \quad I'_f &=
\ideal{t-uf, \p_1+u\textstyle\frac{\p f}{\p x_1}\p_t, \ldots,
  \p_n+u\textstyle\frac{\p f}{\p x_n}\p_t, uv-1} \subseteq D_{x,t}[u,v].
\end{split}
\end{align}
Now outright eliminate $u,v$.
Note that $I'_f$ is quasi-homogeneous if the weights are
$t,u\leadsto -1$ and
$\p_t,v\leadsto 1$, all other variables having weight zero. The
homogeneity of the input and the relation $[\p_t,t]=1$ assures the termination of the computation. The operators of weight $0$ in the output (with  $-\del_tt$ replaced by $s$) generate $I'_f \cap D[s]$.

A modification given in
\cite{Briancon-Maisonobe:BernsteinIdeal} and used, \emph{e.g.},
in~\cite{CastroUcha-JSC04},
reduces the number of algebra generators by one. Consider the
subalgebra $D\<s,\p_t\> \subset D_{x,t}$; the relation
$[s,\p_t]=\p_t$ shows that it
is of solvable type.
According to~\cite{Briancon-Maisonobe:BernsteinIdeal},
\begin{align}
  \begin{split}
\ann_{D[s]}(f^s) &= I''_f \cap D[s],\\ \quad I''_f &=
\ideal{s+f\p_t, \p_1+\textstyle\frac{\p f}{\p x_1}\p_t, \ldots,
  \p_n+\textstyle\frac{\p f}{\p x_n}\p_t} \subset D\<s,\p_t\>.
  \end{split}
\end{align}
Note that $I''_f = \ann_{D_{x,t}}(f^s) \cap D\<s,\p_t\>$. The
elimination step is done as in \cite{Oaku-local-b}; the decrease of
variables usually improves performance.
An algorithm to decide $(A_1)$ for arrangements is given in
\cite{AlvarezCastroUcha-JA07}.

\subsection{Algorithms for the \BS\ polynomial}\label{subsection:BSalgo}
As the minimal polynomial of $s$ on $\calN_f(s)$, $b_f(s)$ can be
obtained by means of linear algebra as a syzygy for the normal forms
of powers of $s$ modulo $\ann_{D[s]}(f^s)+D[s]\cdot f$ with respect to
any fixed monomial order on $D[s]$. Most methods follow this path,
starting with \cite{Oaku-local-b}. Variations appear in
\cite{Walther-JPAA99,OakuTakayama-JPAA01,OTW-JSC00}; see also
\cite{HyperGeomBook}.

A slightly different approach is to compute $b_f(s)$ without recourse
to $\ann_{D[s]}(f^s)$, via a \Grob\ deformation of the ideal
$I_f=\ann_{D_{x,t}}(f^s)$ in~(\ref{eq:ann-Dxt})
with respect to the weight $(-w,w)$ with $w=(0^n,1)\in\NN^{n+1}$:
$\ideal{b_f(s)} = \ini_{(-w,w)}(I_f) \cap \QQ[-\p_tt]$.  Here again,
computing the minimal polynomial using linear algebra tends to provide
some savings in practice.

In \cite{LevandovskyyMorales-JA12} the authors give a method to check
specific numbers for being in $\rho_f$. A method for $b_f(s)$ in the
prehomogeneous vector space setup is in \cite{Muro-Japanese00}.

\subsection{Stratification from $b_f(s)$}
The \Grob\ deformation $\ini_{(-w,w)}(I_f)$
in~\S\ref{subsection:BSalgo} can be refined as follows, see
\cite[Thm.~2.2]{Berkesch-Leykin:MultiplierIdeals}. Let $b(x,s)$
be nonzero in the polynomial ring $\CC[x,s]$.  Then $b(x, s) \in
(\ini_{(-w,w)} I_f) \cap \CC[x,s]$ if and only if there exists $P\in D[s]$
satisfying the functional equation $b(x,s)f^s=Pff^{s}$.  From
this one can design an algorithm not only for computing the {\em
  local} \BS\ polynomial $b_{f,p}(s)$ for $p\in\Var(f)$, but also the
stratification of $\CC^n$ according to local \BS\ polynomials;
see~\cite{Nishiyama-Noro:stratification-by-local-b,
  Berkesch-Leykin:MultiplierIdeals} for various approaches.
Moreover, one can compute the stratifications from Subsection
\ref{subsec-strat-bfu},  see \cite{Leykin-JSC01}.

For the ideal case, \cite{AndresLevandovskyyMorales-Proc09} gives a
method to compute an intersection of a left ideal of an associative
algebra over a field with a subalgebra, generated by a single
element. An application is a method for the computation of the
Bernstein-Sato polynomial of an ideal. Another such was given by
Bahloul in \cite{Bahloul-JSC01}, and a version on general varieties in
\cite{Bahloul-Japan03}.

\subsection{Multiplier ideals}
Consider polynomials $f_1,\dots,f_r \in \CC[x]$, let $f$ stand for
$(f_1,\ldots,f_r)$, $s$ for $s_1,\ldots,s_r$, and $f^s$ for
$\prod_{i=1}^r f_i^{s_i}$.  In this subsection, let $D_{x,t} =
\CC\<x,t,\p_x,\p_t\>$ be the $(n+r)$-th Weyl algebra.

Consider  $D_{x,t}(s)\bullet f^s\subseteq
R_{x,t}[f^{-1},s]f^s$ and put
\begin{align*}
  t_j \bullet h(x, s_1, \ldots, s_j, \ldots , s_r) f^s &= h(x, s_1,
  \ldots, s_j+1, \ldots , s_r) f_j f^s,\\ \p_{t_j} \bullet h(x, s_1,
  \ldots, s_j, \ldots , s_r) f^s &= -s_j h(x, s_1, \ldots, s_j-1,
  \ldots , s_r) f_j^{-1} f^s,
\end{align*}
for $h \in \CC[x][f^{-1},s]$, generalizing the univariate constructions.

The \emph{generalized \BS\ polynomial} $b_{f,g}(\sigma)$ of $f$ at
$g\in\CC[x]$ is the monic univariate polynomial $b$ of the lowest
degree for which there exist $P_k \in D_{x,t}$ such that
\begin{align}\label{eq:generalB}
b(\sigma)gf^s = \sum_{k=1}^r P_k g f_k f^s,\quad\sigma = -\left(
\sum_{i=1}^r \p_{t_i}t_i \right).
\end{align}

An algorithm for $b_{f,g}(\sigma)$ is an essential
ingredient for the algorithms
in~\cite{Shibuta:multiplier-ideals,Berkesch-Leykin:MultiplierIdeals}
that compute the jumping numbers and corresponding multiplier ideals
for $I = \ideal{f_1,\ldots,f_r}$. That $b_{f,g}(\sigma)$ is related to
multiplier ideals was worked out in \cite{BudurMustataSaito-Compos06}.

There are algorithms for special cases: monomial
ideals~\cite{Howald-TAMS01}, hyperplane
arrangements~\cite{Mustata-TAMS06}, and determinantal
ideals~\cite{JohnsonA.-thesis}. A {\em Macaulay2} package {\em
  MultiplierIdeals} by Teitler collects all available (in {\em
  Macaulay2}) implementations. See also \cite{Budur-Vfilt05}.

\subsection{Software}
Algorithms for computing \BS\ polynomials have been implemented in
{\em kan/sm1}~\cite{KANwww}, {\em Risa/Asir}~\cite{risa-asir-www},
{\em dmod\_lib} library \cite{Dmods-Singular}
for {\em Singular}~\cite{DGPS}, and the {\em
  D-modules} package~\cite{DmodulesM2} for {\em
  Macaulay2}~\cite{M2www}.  The best source of information of these is
documentation in the current versions of the corresponding software.
A relatively recent comparison of the performance for several families
of examples is given in~\cite{Levandovskyy-Morales:comparison}.

The following are articles by
developers discussing their implementations:
\cite{NoroBPoly,Nishiyama-Noro:stratification-by-local-b,OakuTakayama-JPAA01,
  Singular:constructive-D-modules,Dmods-Singular, Leykin2002Dmodules,
  Berkesch-Leykin:MultiplierIdeals}.


\bibliographystyle{plain}
\bibliography{bib}

\end{document}